\documentclass[12pt]{article}

\usepackage{epsfig}

\author{Manjunath Krishnapur and Yuval Peres\footnote{Research supported
in part by NSF grants
\#DMS-0104073 and \#DMS-0244479. Part of this work was done while the
second author was visiting Microsoft Research.}}
\title{Recurrent graphs where two independent random walks
\\ collide finitely often}

\newtheorem{theorem}{Theorem}[section]
\newtheorem{lemma}[theorem]{Lemma}

\newenvironment{proof}[1][Proof]{\begin{trivlist}
       \item[\hskip \labelsep {\bfseries #1}]}{\end{trivlist}}
       \newenvironment{definition}[1][Definition]{\begin{trivlist}
       \item[\hskip \labelsep {\bfseries #1}]}{\end{trivlist}}
       
       \newenvironment{remark}[1][Remark]{\begin{trivlist}
       \item[\hskip \labelsep {\bfseries #1}]}{\end{trivlist}}
\newcommand{\qed}{\nobreak \ifvmode \relax \else
             \ifdim\lastskip<1.5em \hskip-\lastskip
             \hskip1.5em plus0em minus0.5em \fi \nobreak
             \vrule height0.75em width0.5em depth0.25em\fi}

\def\wed{\wedge}
\def\E{{\bf E}}
\def\P{{\bf P}}
\def\v{{\bf v}}

\def\Gam{\Gamma}

\def\alp{\alpha}

\def\eps{\epsilon}

\def\given{\arrowvert}

\def\summ{\sum\limits}

\def\l{\left}
\def\r{\right}

\def\hsp{\hspace}
\def\mb{\mbox}
\def\G{{\bf{G}}}
\def\HH{{\bf{H}}}
\def\Z{{\bf{Z}}}

\begin{document}
\maketitle
\begin{abstract} We present a class of graphs where simple random walk
  is recurrent, yet two independent walkers meet only finitely many
  times almost surely. In particular, the comb lattice, obtained from $\Z^2$ by
  removing all horizontal edges off the $x$-axis, has this property.  We also conjecture
  that the same property holds
  for some other graphs, including the incipient infinite cluster for
  critical percolation in $\Z^2$.
\end{abstract}

\section{Introduction}
In ``Two Incidents''~\cite{polya}, George P\'{o}lya describes the
incident that led him to his celebrated results on random walks on Euclidean
lattices:
\begin{quote} \it
``$\ldots$ he and his fianc\'{e}e (would) also set out for a stroll in the
woods, and then suddenly I met them there. And then I met them the
same morning repeatedly, I don't remember how many times, but
certainly much too often and I felt embarrassed: It looked as if I was
snooping around which was, I assure you, not the case.
I met them by accident - but how likely was it that it happened by
accident and not on purpose?''
\end{quote}
P\'{o}lya formulated the problem of the meeting of two walkers for
 random walks on a Euclidean lattice; in that case, it
reduces to the problem of a single walker returning to his starting point.
As we show in this paper, these two problems can have different answers when the
 ambient graph is not transitive.

 Call a graph $\G$ recurrent if simple random walk on it is
  recurrent. Say that a graph $\G$ has the {\bf finite
  collision property} if two independent simple random walks $X,Y$ on
  $\G$ starting from the same vertex meet only finitely many times,
  i.e., $|\{n: X_n=Y_n\}|<\infty$, almost surely. Our goal
 is to present a class of recurrent graphs with the
finite collision property.

\begin{figure}
\begin{center}
\epsfig{file=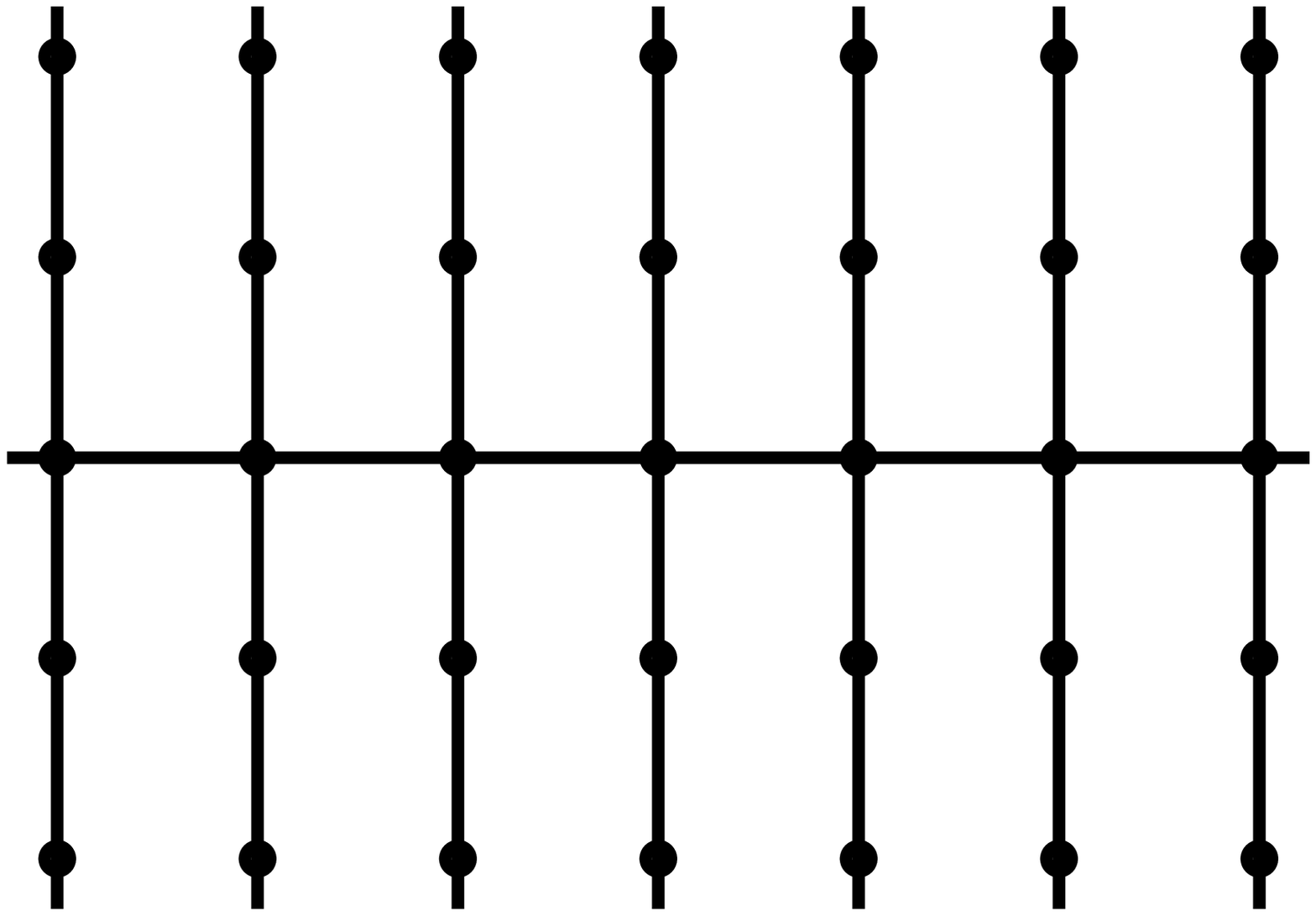, height=7cm, width=8cm}
\end{center}
\caption{Comb lattice ($\mbox{Comb}(\Z)$) }
\label{fig:comblattice}
\end{figure}

\begin{definition} Given a graph $\G$, let $\mbox{Comb}(\G)$ be the graph with vertex set $V(\G)\times \Z$ and edge set
\[ \l\{\l[(x,n),(x,m)\r]: |m-n|=1
  \r\}\cup\l\{\l[(x,0),(y,0)\r]:\l[x,y\r] \mbox{ is an edge in
  $\G$}\r\}. \]
\end{definition}
In words, this means that we attach a copy of $\Z$ at each vertex of
  the graph $\G$. See Fig~\ref{fig:comblattice} for a picture of $\mbox{Comb}(\Z)$.
Clearly, if $\G$ is recurrent, so is $\mb{Comb}(\G)$.

\begin{theorem}\label{thm:zcombs} Let $\G$ be any recurrent infinite
  graph with constant vertex degree. Then  $\mbox{Comb}(\G)$  has the finite collision property.
\end{theorem}

We note the following points.
\begin{itemize}
\item Liggett~\cite{liggett} has given examples of symmetric recurrent Markov chains
for which two independent copies of the chain collide only finitely many times. Those examples
are not simple random walks  on graphs, however.

\item If $X,Y$ are independent random walks on a graph starting from a vertex $v$, then the expected number of meetings between them is $\summ_{n}
 \summ_{w}(p^{(n)}(v,w))^2$, where $p^{(n)}$ is the $n$-step transition
 function. Now, \[ p^{(2n)}(v,v)=\summ_w (p^{(n)}(v,w))^2{\pi(v)\over
 \pi(w)}, \] where $\pi(w)$ denotes the degree of $w$. Therefore, if the expected number of meetings between two independent walkers is finite, so is $\summ_n p^{(2n)}(v,v)$, and hence the graph is transient.

The converse is true for bounded degree graphs because then ${\pi(v)
  \over \pi(w)}$ is bounded away from zero. This converse can fail
  when the degrees are unbounded as the following example shows:
  Consider $\Z_+$ and add $2^n$ disjoint paths of length $2$ between
  $n$ and $n+1$. Then, $X_{2n}$ and $Y_{2n}$ are just random walks on $\Z_+$
  with a bias of ${1\over 3}$ to the right and hence the graph is transient. However they meet infinitely often almost surely (the difference eventually coincides with an unbiased random walk on $\Z$ and hence visits zero infinitely often).

\item Recurrent transitive graphs cannot have the finite collision
  property. This is because transitivity clearly implies
  that the number of meetings has a Geometric distribution, whence it
  is finite only if it has finite expectation.
\end{itemize}

\section{Proof of Theorem~\ref{thm:zcombs}}

{\bf Proof of Theorem~\ref{thm:zcombs} }
As noted earlier, $\mb{Comb}(\G)$ is recurrent. We only need to prove the
finite collision property. Let $X$ and $Y$ be independent simple random walks
(SRWs) on $\mbox{Comb}(\G)$ starting from the same vertex $(o,0)$. We  make the
following definitions.
 \[ \begin{array}{l}\label{arr:defn1}
 Z_{n,\ell} \colon= |\{(N,L): n\le N\le 2n, \ell\le L \le 2\ell \mbox{ and }
X_N=Y_N=(v,L) \mbox{ for some }v\in\G \}|.\\
 A_{n,\ell}\colon= \{ Z_{n,\ell}>0\}.\\
 W_{n,\ell}\colon=\sum\limits_{k={\ell \over 2},\ell,2\ell}Z_{n,k}+\sum\limits_{k={\ell \over 2},\ell,2\ell}Z_{2n,k}. \\
 \end{array} \]
(Here $|S|$ denotes the number of elements of $S$.)

Let $d$ denote the common degree of vertices in $\G$. In what follows, $C,C_1,C_2$ etc. will denote positive finite constants whose values
may change from one appearance to another.

\begin{lemma}\label{lem:lemma1} $\E[Z_{n,\ell}]\le C \ell n^{-1/4}$ for some
  finite constant $C$, $\forall n,\ell\ge 1$.
\end{lemma}

\begin{remark} For the case when $\G=\Z$, i.e., for $\mbox{Comb}(\Z)$, this
  lemma is suggested by the fact that
\begin{equation}\label{eq:woess}
p^{(2n)}(0,0)\sim {\sqrt{2}\over\Gam(1/4)}n^{-3/4}.
\end{equation}
See Proposition 18.4 in Woess~\cite{woess} for a proof of (\ref{eq:woess}).
\end{remark}

\begin{proof} (Lemma~\ref{lem:lemma1}) We generate the random walks $X$ and $Y$ in the
  following manner. Let $U$ and $U'$ be independent simple random walks on
  $\G$ starting from a vertex $o$. Let $V$ and $V'$ be independent simple
  random walks on $\Z$, with the modification that they have a
  self-loop probability of ${d \over d+2}$ at
  $0$. Let $K_n$ and  $K'_n$ be the number
  of transitions of $V$ and
  $V'$ from $0$ to $0$ in the first $n$ steps. Then set
\[ X_n=(U_{K_n},V_n)  \mb{ and } Y_n=(U'_{K'_n},V'_n). \]
 It is clear that $X$ and $Y$ are independent simple random walks on
 $\mb{Comb}(\G)$, both starting from $(o,0)$.

 Now fix any $L\in \Z$ and consider
 \begin{eqnarray*}
   \P[X_n=Y_n=(v,L) \mb{ for some } v\in \G] &=& \summ_{k=0}^{\infty}
   \summ_{k'=0}^{\infty} \P[V_n=V'_n=L; K_n=k, K'_n=k'; U_k=U'_{k'}] \\
&=& \summ_{k,k'} \P[V_n=V'_n=L;K_n=k,K'_n=k'] \P[U_k=U'_{k'}]. \\
 \end{eqnarray*}

Given two paths ${\cal P}_1,{\cal P}_2$ of lengths $i$
and $j$ in $\G$ starting from $o$ and having the same endpoint $w$, let ${\cal P}$ be the path
obtained by traversing ${\cal P}_1$ first and returning to $o$ via ${\cal
  P}_2$. Then
\[ \P[\{U_k\}_{k\le i+j}={\cal P} ] = \P[\{U_k\}_{k\le
 i}={\cal P}_1] \P[\{U'_k\}_{k\le j}={\cal P}_2]. \]
by our assumption of constant degrees.
Summing over all possible $w$ and all ${\cal P}_1,{\cal P}_2$, we get
\begin{equation}
\P[U_i=U'_j\given U_0=U'_0=v] = \P[U_{i+j}=v\given U_0=v].
\label{eq:looparound}
\end{equation}

Moreover, for SRW $U$ on any infinite graph with bounded degrees,
 \begin{equation}
  \P[U_n=o|U_0=o]\le {C \over \sqrt{n}},
 \label{eq:rootnbound}
 \end{equation}
for some constant $C$ (not depending on $o$). For a proof, see
Woess~\cite{woess}, Corollary 14.6 .

From (\ref{eq:looparound}) and (\ref{eq:rootnbound}) we get
\begin{eqnarray}
  \label{eq:ubd}
 \P[X_n=Y_n=(v,L) \mb{ for some } v\in \G] & \le & 
  \P[V_n=V'_n=L;K_n=K'_n=0] \\[1ex] &+& C \E\l[ {{\mathbf 1}(V_n=V'_n=L;K_n+K'_n>0)\over
  \sqrt{K_n+K'_n}} \r]. \nonumber
\end{eqnarray}

To bound this quantity we think of $V_n$ as being generated in the
following manner. Take a simple random walk $\{S_n\}$ on $\Z$ (no
self-loop at $0$) starting from $0$ and let $\{G_i\}$ be i.i.d.
Geometric(${d\over d+2}$) random variables. To be precise, this means that
$\P[G_i=k] = ({d \over d+2})^k{2\over d+2}$ for $k\ge 0$.
 Then we generate $V$ by
following the path $S$ except that at the $i$th visit to the origin by $S$,
the walk $V$ stays there for $G_i$ steps before taking the next step
according to $S$. Similarly, $V'$ is generated using $S'$ and
$\{G'_i\}$.

Then let $H_n=\summ_{i=1}^{n/2} {\mathbf 1}(S_i=0)$ and similarly define
$H'_n$. Then, either $K_n\ge R_n:=\summ_{i=1}^{H_n} G_i$ or else
$K_n \ge {n\over 2}$. Therefore the second summand on the right in
(\ref{eq:ubd}) can be bounded by (omitting the constant $C$)
\begin{equation}
   \E\l[{{\mathbf 1}(V_n=V'_n=L;R_n+R'_n>0)\over \sqrt{(R_n\wedge {n\over 2})+(R'_n\wedge {n\over 2})}}\r] + \P[V_n=V'_n=L;R_n=R'_n=0].
\label{eq:ubd2}
\end{equation}
%\\ \qquad{}+ \P[\max\{R_n,R'_n\}\ge {n\over 2};V_n=V'_n=L]
Condition on $\{S_i:i\le n/2\}$ and $\{G_i:i\le H_n\}$ and on
their primed counterparts. If it happens that $\max\{R_n,R'_n\}<{n\over 4}$,
then $V$ and $V'$ have at least $n/4$ more steps to go and hence the
conditional probability that $V_n=V'_n=L$ is at most ${C^2\over n}$
(because $\P[V_{n/4}=L']\le {C\over \sqrt{n}}$ for any $L'$).
Thus the first term in (\ref{eq:ubd2}) can be bounded by
\begin{equation}
 C' \E\l[{{\mathbf 1}(R_n+R'_n>0)\over \sqrt{R_n+R'_n}}{1\over n} + {{\mathbf 1}\l(\max\{R_n,R'_n\}>{n\over 4}\r)\over \sqrt{n}}  \r].
\label{eq:ubd3}
\end{equation}

We recall the following facts
\begin{itemize}
\item $\E\l[H_n^{-1/2}{\mathbf 1}\l(H_n\ge 1\r)\r]\le C_1 n^{-1/4}$. To see this, consider
\begin{eqnarray*}
\E\l[H_n^{-1/2}{\mathbf 1}\l(H_n\ge 1\r)\r] &\le& {1\over n^{1/4}}+\summ_{k=1}^{n^{1/2}}{\P[H_n=k]\over \sqrt{k}} \\
 &\le& {1\over n^{1/4}}+{C_1\over n^{1/2}}\summ_{k=1}^{n^{1/2}} {1\over \sqrt{k}} \hspace{1.5cm} \l(\mb{ since }\P[H_n=k]\le {C\over n^{1/2}}\forall k\r)\\
 &\le& C_1 n^{-1/4}.
\end{eqnarray*}
\item If $\{G_i\}$ are i.i.d. Geometric($p$) random variables, then
  \begin{equation}\label{eq:geometrics}
 \E\l[{{\mathbf 1}\l(\summ_{i=1}^r G_i \not=0 \r)\over \sqrt{\summ_{i=1}^r
      G_i}}\r]\le {C(p)\over \sqrt{r}}.
 \end{equation}
(Let $\mu=\E[G_i]$. If ${1\over r}\summ_{i=1}^r G_i >{\mu - \eps}$, then the random variable in (\ref{eq:geometrics}) is
less than ${1\over \sqrt{r(\mu - \eps)}}$. The probability that ${1\over r}\summ_{i=1}^r G_i$ is less than ${\mu - \eps}$ decays exponentially, by Cram\'er's theorem).
\end{itemize}
These facts immediately give
\[ \E\l[{{\mathbf 1}(R_n+R'_n>0)\over \sqrt{R_n+R'_n}}\r] \le {C_3 \over n^{1/4}}. \]
We ultimately want to get a bound for $\P[X_n=Y_n=(v,L) \mb{ for some
} v\in \G]$. From what we have done so far this is bounded by the sum
of the following three terms
\begin{itemize}
\item The first term in (\ref{eq:ubd3}) is bounded by $C_4 n^{-5/4}$.
\item The second term in (\ref{eq:ubd3})
is bounded by $C_5\P[R_n>{n\over 4}]/\sqrt{n}$. This decays super-polynomially.
\item The second term in (\ref{eq:ubd2}) and the first term in
  (\ref{eq:ubd}) are together bounded by \[C_6\P[V_n=V'_n=L,R_n=R'_n=0].\]
 To bound $\P[V_n=V'_n=L,R_n=R'_n=0]$,
condition on $\{S_i:i\le {n\over 2}\}$, $\{G_i:i\le {n\over 2}\}$ and
their primed versions as before. Since  the probability that a simple
random walk on $\Z$ does not return to  zero up to time $n$ is
asymptotic to ${C\over \sqrt{n}}$, we can easily deduce that
\begin{eqnarray*}
\P[V_n=V'_n=L \mb{ and }R_n=R'_n=0] = O\l({1\over n^2}\r).
\end{eqnarray*}
\end{itemize}

Thus we get
\begin{equation}
  \P[X_n=Y_n=(v,L) \mb{ for some } v\in \G] \le {C \over n^{5/4}}
  \hsp{1.5cm} \mb{ for every } L\in \Z, n\ge 1.
\label{eq:ubdfortransprob}
\end{equation}

Now,\begin{eqnarray*}
\E[Z_{n,\ell}] &=& \summ_{N=n}^{2n}
\summ_{L=\ell}^{2\ell}\P[X_N=Y_N=(v,L) \mbox{ for some } v\in\G] \\
&\le& n\ell {C \over n^{5/4}} = C {\ell \over n^{1/4}},
\end{eqnarray*}
as claimed. $\qed$
\end{proof}

Note that from the above lemma we also get
\begin{equation}
\E[W_{n,\ell}]\le C\frac{\ell}{n^{1/4}} \mbox{ \hspace{1cm} for all
 $n,\ell \ge 1$ and some constant } C<\infty.
\label{eq:upperbound}
\end{equation}

\begin{lemma}\label{lem:lemma2} Fix $0<\alp<1$. There is a constant $C>0$
(depending on $\alp$ but not on $\ell$ or $n$) such that for
all $n,\ell$ with $1\le \ell<2(2n)^{1/2\alp}$, we have $\E[W_{n,\ell}|
A_{n,\ell}]\ge C\ell^{\alp}$.
\end{lemma}

\begin{proof} (Lemma~\ref{lem:lemma2}) Suppose $A_{n,\ell}$ occurs. Then the two random
walks collide at a time $N$ with $n\le N\le 2n$, and at some vertex $(v,L)$ with
$\ell \le L \le 2\ell$. Then in $W_{n,\ell}$ we are counting all
collisions that occur for the next $2n$ steps or till one of the walks
reaches $(v,L\pm{\ell \over 2})$, whichever occurs earlier.
 By considering only collisions that occur before one of them hits
 $(v,L\pm{\ell \over 2})$ the problem is reduced to one about random walks on
 a segment of $\Z$.

More precisely, let $U,V$ be two independent random
 walks on $\Z$ starting from $0$. Let $T_U$ be the first time $U$
 hits $\pm {\ell\over 2}$ and similarly define $T_V$. If
\[ Y_{n,\ell}=\sum\limits_{k=0}^{2n\wed T_U \wed T_V }
 {\mathbf 1}(U_k=V_k), \]
then given that the event $A_{n,\ell}$ occurs, $W_{n,\ell}$ is stochastically larger than $Y_{n,\ell}$. Therefore, if $2n\ge (\ell /2)^{2\alp}$, then
\begin{eqnarray*}
\E[W_{n,\ell}\given A_{n,\ell}] &\ge& \E[Y_{n,\ell}]\\
&\ge& \sum\limits_{k=0}^{(\ell/2)^{2\alp}} \P[U_k=V_k;
T_U\wed T_V>k] \\
   &\ge& \l( \sum\limits_{k=0}^{(\ell/2)^{2\alp}} \P[U_k=V_k] \r) -
   (\ell/2)^{2\alp}\P[T_U\wed T_V\le(\ell/2)^{2\alp}] \\
  &\ge& \sum\limits_{k=0}^{(\ell/2)^{2\alp}} \frac{C'}{\sqrt{k}} -
  (\ell/2)^{2\alp}2\P[T_U\le(\ell/2)^{2\alp}],
\end{eqnarray*}
since for independent SRWs $U,V$ on $\Z$, we have $\P[U_k=V_k]\sim C' k^{-1/2}$. Observe that $\P[T_U\le(\ell/2)^{2\alp}]$ tends to zero faster than any
  polynomial in $\ell$. Therefore,
\[ \E[W_{n,\ell}\given A_{n,\ell}] \ge C\ell^{\alp}. \]

This proves the lemma. $\qed$
\end{proof}

From the two lemmas above, given $\alp<1$ we have constants $C_1,C_2$ such that
\begin{eqnarray}
  \E[W_{n,\ell}] &\le& C_1\frac{\ell}{n^{1/4}} \mbox{\hspace{1cm} for every
    $\ell,n$} \label{eq:claim1}, \\
  \E[W_{n,\ell}|A_{n,\ell}] &\ge& C_2\ell^{\alp} \mbox{\hspace{1cm} for $\ell\le
  2(2n)^{1/2\alp}$} \label{eq:claim2},
\end{eqnarray}
whence we get
\begin{equation}\label{eq:ubdforA}
\P[A_{n,\ell}] \le \frac{\E[W_{n,\ell}]}{\E[W_{n,\ell}|A_{n,\ell}]}
\le C\frac{\ell^{1-\alp}}{n^{1/4}} \mbox{\hspace{1.5cm} for $\ell\le 2(2n)^{1/2\alp}$},
\end{equation}
for yet another constant $C$.

 Now we let $n,\ell$ satisfying $\ell\le 2(2n)^{1/2\alp}$  run over powers
 of $2$, and get
\begin{eqnarray*}
\sum\limits_{r=0}^{\infty}\sum\limits_{k=0}^{1+\frac{r+1}{2\alp}}
\P[A_{2^r,2^k}] &\le&  \sum\limits_{r=0}^{\infty}\sum\limits_{k=0}^{1+\frac{r+1}{2\alp}} C\frac{2^{k(1-\alp)}}{2^{r/4}} \hspace{2cm} \mbox{ by
  (\ref{eq:ubdforA}) } \\
  &\le& \sum\limits_{r=0}^{\infty} C' \frac{2^{r(1-\alp)/2\alp}}{2^{r/4}} \\
    &<& \infty \mbox{\hspace{2cm} if $\alp > \frac{2}{3}$}.
\end{eqnarray*}
Thus almost surely only finitely many of the events $A_{2^r,2^k}$ for $k\le
 1+\frac{r+1}{2\alp}$ occur. This shows that if $2/3< \alp <1$, then the set
\[ \{n: X_n=Y_n=(v,\ell) \mbox{\hspace{.1cm} for
some }v\in \G \mbox{ and $\ell$ with } 1\le |\ell|\le 2(2n)^{1/2\alp} \} \] is finite almost surely (since each such $(n,\ell)$
is contained in one of the above sets). We proved this for $1\le \ell \le 2(2n)^{1/2\alp}$. By symmetry, the same holds for negative $\ell$. The number of meetings on the backbone, i.e., on $\{\ell=0\}$, is finite, by (\ref{eq:ubdfortransprob}).

However, $\{n: |V_n|>2(2n)^{1/2\alp} \mbox{ or }
|V'_n|>2(2n)^{1/2\alp}\}$ is finite almost surely as can be easily seen,
for instance, from the law of iterated logarithm.
\\

This proves that the total number of collisions between the two random
walkers on $\mb{Comb}(\G)$ is finite almost surely. $\qed$

\section{More Examples}
\begin{definition} Given two graphs $\G ,\HH$, and a vertex $\v$ of $\HH$,
  define $\mb{Comb}_{\v}(\G,\HH)$ to be the
  graph with vertex set $V(\G)\times V(\HH)$ and edge set
\[ \l\{\l[(x,w),(x,z)\r]:\l[w,z\r] \mbox{ is an edge in \HH}
  \r\}\cup\l\{\l[(x,\v),(y,\v)\r]:\l[x,y\r] \mbox{ is an edge in
  $\G$}\r\}. \]
When $\HH=\Z$ (and without loss of generality $\v=0$), this clearly reduces to
$\mb{Comb}(\G)$.
\end{definition}

If $\G,\HH$ are recurrent, and $\v$ is a vertex of $\HH$, then $\mb{Comb}_{\v}(\G,\HH)$ is also obviously recurrent.
When $\HH=\Z^2$, we take $\v=(0,0)$ and drop the subscript $\v$ in $\mb{Comb}_{\v}(\G,\HH)$.

\begin{theorem} \label{thm:z2combs} Let $\G$ be any recurrent infinite
  graph with constant vertex degree. Then  $\mbox{Comb}(\G,\Z^2)$ has the finite collision property.
\end{theorem}

\begin{proof} As the proof is very similar to that of
  Theorem~\ref{thm:zcombs} (the difference is in the estimates) we
  shall only briefly sketch the main steps.

For $\ell\ge 1$, let $B_\ell=\{(x,y)\in \Z^2: \ell\le \max\{|x|,|y|\} \le 2\ell \}$ be
the annulus of radii $\ell$ and $2\ell$. Then we define
\[ Z_{n,\ell} =   |\{(N,L): n\le N\le 2n, L \in B_{\ell} \mbox{ and }
X_N=Y_N=(v,L) \mbox{ for some }v\in\G \}|. \]
Then define $A_{n,\ell}$ and $W_{n,\ell}$ as before.
 Then analogously to Lemma~\ref{lem:lemma1} and
Lemma~\ref{lem:lemma2} we have the following lemma.
\begin{lemma}\label{lem:lemma3} With the above definitions,
\begin{itemize}
\item $\E\l[W_{n,\ell}\r] \le C_1{\ell^2 \over n\sqrt{\log(n)}}$
  $\forall 1\le \ell,n$.
\item $\E\l[W_{n,\ell} \given A_{n,\ell} \r]\ge C_2\log(\ell)$ for
  $1\le \ell \le n$.
\end{itemize}
\end{lemma}

\begin{proof} (Lemma~\ref{lem:lemma3}) The upper bound for $\E\l[W_{n,\ell}\r]$ can be proved along the same lines as Lemma~\ref{lem:lemma1}. First we prove
\begin{equation} \label{eq:ubd2forz2}
  \P\l[X_n=Y_n=(v,L) \mb{ for some } v\in \G \r] \le {C \over
  n^2\sqrt{\log(n)}}.
\end{equation}
%(The reason is simple: In $n$ steps of the walk on $\G\oplus \Z^2$,
%the number of ``vertical'' steps is of the order of $n$, and since SRW on
%$\Z^2$ returns an order of $\log(n)$ times to the starting point in
%$n$ steps, the number of ``horizontal'' steps is also of the same
%order.)

All the steps go through without change till (\ref{eq:ubd3}). Moreover, the terms with
$R_n=0$ or $R_n\ge {n\over 4}$ etc can be shown to be of lower order in
the same manner. (To bound the terms with $\{R_n=0\}$, use the fact that for simple random walk in the plane, $\P[H_n=0]\le {C\over \log n}$. See Erd\H{o}s and Taylor~\cite{et}.) Only note that in $\Z^2$ the $n$-step transition
probabilities are bounded by $Cn^{-1}$. The dominant term is
\[ \E\l[{{\mathbf 1}(V_n=V'_n=L;R_n+R'_n>0)\over \sqrt{(R_n\wedge {n\over 2})+(R'_n\wedge {n\over 2})}}\r], \]
where the notations are as before (now $V,V'$ are random walks on
$\Z^2$ instead of $\Z$).

For simple random walk in the plane $\P[H_n=k]\le {C\over \log n}$
$\forall n,k$ ($H_n$ is the number of returns to origin by time $n$.
 See Erd\H{o}s and Taylor~\cite{et}).
Using this and the bound (\ref{eq:geometrics}) for i.i.d. Geometric variables,  we get the bound (\ref{eq:ubd2forz2}). In $W_{n,\ell}$ we are
  counting (up to constants) $n$ steps and $\ell^2$ sites, and thus the upper bound for
  $\E\l[W_{n,\ell}\r]$ follows.

The lower bound for $\E\l[W_{n,\ell}\given
A_{n,\ell}\r]$ is even easier. Referring back to the proof of
Lemma~\ref{lem:lemma2}, a lower bound can be obtained by counting only
those meetings that occur for a duration of $2n$ and before one of the
two walkers goes a distance of $\ell/2$ from the meeting point (that is
assured by $A_{n,\ell}$). Since at least $\ell/2$ steps are needed to
go a distance $\ell/2$, if $4n>\ell$,
\begin{eqnarray*}
 \E\l[W_{n,\ell}\given A_{n,\ell}\r] &\ge&
 \summ_{k=1}^{\ell/2}\P[U_k=V_k] \hsp{1cm} U,V \mb{ are SRWs on }\Z^2\\
 &\ge& \summ_{k=1}^{\ell/2}{C'\over k} \ge C_2\log(\ell).
\end{eqnarray*}
This completes the proof of Lemma~\ref{lem:lemma3}. $\qed$
\end{proof}

From Lemma~\ref{lem:lemma3} we get
\begin{eqnarray*}
  \P\l[A_{n,\ell}\r] &\le&  {\E[W_{n,\ell}] \over
    \E[W_{n,\ell}\given A_{n,\ell}]} \\
    &\le& C{\ell^2 \over n\sqrt{\log(n)}\log(\ell)} \hsp{2cm} \mb{ for
    }2\le \ell<4n.
\end{eqnarray*}

Now we let $n,\ell$ run over powers of $2$ but only over pairs for
which $2\le \ell \le \sqrt{n} (\log(n))^{1/8}$ (trivially the above bound
 for $\P[A_{n,\ell}]$ holds for these values of $n,\ell$). Here $\log$ denotes logarithm to base $2$. Then
\begin{eqnarray*}
\summ_{r=1}^{\infty}\sum\limits_{k=1}^{{r\over 2}+{1\over 8}\log(r)}
\P[A_{2^r,2^k}] &\le&
\summ_{r=1}^{\infty}\sum\limits_{k=1}^{{r\over 2}+{1\over 8}\log(r)}
C{2^{2k} \over 2^r \sqrt{r} k}  \\
  &\le& \summ_{r=1}^{\infty} {C_0 \over r^{5/4}}\\
    &<& \infty,
\end{eqnarray*}
where, in the penultimate line, we have used the following easily checked
fact:
\[ \summ_{k=1}^n {4^k \over k} \le C{4^n \over n} \]
for some constant $C$ not depending on $n$.

This proves that almost surely only finitely many of the events
$A_{2^r,2^k}$ with $k\le {r\over 2}+{1\over 8}\log(r)$ occur.(The
cases $\ell=0,1$ are taken care of directly by (\ref{eq:ubd2forz2}).)
However, as before, letting $V_n=(V_n^{(1)},V_n^{(2)})$ and similarly for $V'$, we observe that $\{n: \max\{|V_n^{(1)}|,|V_n^{(2)}|\}>\sqrt{n}(\log(n))^{1/8} \mbox{ or }
\max\{|V_n^{'(1)}|,|V_n^{'(2)}|\}>\sqrt{n}(\log(n))^{1/8}\}$ is finite almost surely, as shown by  the law of iterated logarithm.
$\qed$
\end{proof}

\section{Questions}

\begin{itemize}
\item Is it true for any two infinite recurrent graphs $\G,\HH$ and any vertex $\v\in \HH$
  that $\mb{Comb}_{\v}(\G,\HH)$ has the finite collision property?
\item If $\HH_n$ is a sequence of finite graphs then the graph
  obtained by attaching $\HH_n$ to the vertex $n$ of $\Z$ gives a
  comb-like structure similar to the examples given in this paper. This leads us to the following questions.
  \begin{itemize}
  \item Do trees in the uniform and minimum spanning forests on $\Z^d$ have
  the finite collision property? For definitions and properties of Uniform and Minimal Spanning forests see Lyons and Peres~\cite{lp}.
  \item  Does a critical Galton-Watson tree conditioned to survive have the same
  property ? (Assume that the offpring distribution has finite variance.)
   This conditioning on an event of zero probability can be made precise easily; see Kesten~\cite{kesten1}.
   \end{itemize}

  The reason for expecting such behavior is that these trees are known
  to be ``one-ended'', meaning that they have the comb-like structure
  described above (although the ``backbone'' extends infinitely in only one direction).
  \item Does the incipient infinite cluster in $\Z^2$ (this is the cluster containing the origin in bond percolation on $\Z^2$ at criticality, conditioned to be infinite) have the finite collision property? It is known that almost surely there is no infinite cluster in $\Z^2$ at criticality. However, the incipient infinite cluster can still be defined. See Kesten~\cite{kesten2}.
\end{itemize}

{\bf Acknowledgement: } We thank Jeffrey Steif and Nina Gantert for encouragement.

\sc \bigskip \noindent Yuval Peres, Departments of
Statistics and Mathematics, U.C.\ Berkeley, CA 94720,
USA. \\{\sf peres@stat.berkeley.edu}, {\sf
stat-www.berkeley.edu/$\sim$peres}

\sc \bigskip \noindent Manjunath Krishnapur, Department of
Statistics, U.C.\ Berkeley, CA 94720,
USA. \\{\sf manju@stat.berkeley.edu}

\end{document}